\documentclass{article}
\usepackage[utf8]{inputenc}

\title{Nonparametric signal detection with small values of type I and type II error probabilities}

\author{Mikhail Ermakov }
\date{October 2020}

\usepackage{srcltx}
\usepackage{graphicx,color,amssymb,latexsym,amsmath,amsfonts}
\usepackage{mathrsfs}

\RequirePackage{amsthm,amsmath,amssymb}
\RequirePackage[numbers]{natbib}

\numberwithin{equation}{section}
\theoremstyle{plain}
\newtheorem{thm}{Theorem}[section]

\newtheorem{remark}{Remark}[section]

\begin{document}
\global\long\def\zb{\boldsymbol{z}}
\global\long\def\ub{\boldsymbol{u}}
\global\long\def\vb{\boldsymbol{v}}
\global\long\def\wb{\boldsymbol{w}}
\global\long\def\tb{\boldsymbol{t}}
\global\long\def\eb{\boldsymbol{e}}
\global\long\def\sb{\boldsymbol{s}}
\global\long\def\0b{\boldsymbol{0}}
\global\long\def\xb{\boldsymbol{x}}
\global\long\def\xib{\boldsymbol{\xi}}
\global\long\def\etab{\boldsymbol{\eta}}
\global\long\def\omegab{\boldsymbol{\omega}}
\global\long\def\yb{\boldsymbol{y}}
\global\long\def\thetab{\boldsymbol{\theta}}

\maketitle
Institute of Problems of Mechanical Engineering RAS, Bolshoy pr., 61, VO, 1991178  St. Petersburg and
St. Petersburg State University, Universitetsky pr., 28, Petrodvoretz, 198504 St. Petersburg, RUSSIA


keywords :
Neymann test,
consistency,
signal detection,
goodness-of-fit tests

\footnote{Research has been supported RFFI Grant  20-01-00273.}

\begin{abstract}We consider problem of signal detection in Gaussian white noise. Test statistics are linear combinations of      squares of estimators of Fourier coefficients
                                     or $\mathbb{L}_2$-norms of kernel estimators.
We point out
 necessary and sufficient conditions when nonparametric sets of alternatives have
   a given rate of exponential decay for type II                                                                                                    error probabilities.
\end{abstract}

\section{Introduction}
In hypothesis testing and confidence estimation we usually work with events having small error probabilities.
In paper we explore probabilities of  such events arising in problems of
signal detection in Gaussian white noise.
 For tests of Neymann type and for tests  generated $\mathbb{L}_2$-norms of kernel estimators, we provide comprehensive description of nonparametric sets of alternatives having given rates of convergence to zero of type II error probabilities.

 If problem of nonparametric hypothesis testing is the problem of testing hypothesis on adequacy of statistical model choice, the alternatives can be arbitrary. The only natural requirement is that functions from set of alternatives does  not belong to hypothesis. Thus we need to explore consistency of nonparametric sets of alternatives without any assumptions.
We should not be limited by assumptions   of their parametric structure or  ownership to
some class of smooth functions  \cite{er90, ing02}. Such an approach has been developed in \cite{dal13, erm21, erm22} and this line is continued in this paper.

In papers \cite{erm21,erm22} we established necessary and sufficient conditions of uniform consistency of nonparametric sets of alternatives for wide spread nonparametric tests: Kolmogorov, Kramer-von Mises, chi-squared with number of cells increasing with growing of sample size, Neymann, Bickel-Rosenblatt. The results were provided  for nonparametric sets of alternatives assigned both in terms of distribution function and density.

The paper provides a comprehensive description of uniformly consistent nonparametric sets of alternatives having a given rate of convergence of type II error probabilities to zero as noise power tends to zero. Such a description is provided, both in terms of the distance method (asymptotics of distance between hypothesis and alternatives) and in terms of  rate of convergence to zero
  $\mathbb{L}_2$ -- norms of signals belonging to sets of alternatives. As mentioned we consider tests of Neymann type \cite{ney} and  tests  generated $\mathbb{L}_2$-norms of kernel estimators  \cite{bic, hor}.

Let we observe random process $Y_\varepsilon(t)$, $t \in [0,1]$, defined stochastic differential equation
\begin{equation}\label{vv1}
dY_\epsilon(t) = S(t)\, dt + \varepsilon\, dw(t), \quad \varepsilon > 0,
\end{equation}
where $S \in \mathbb{L}_2(0,1)$ and $d\,w(t)$ - Gaussian white noise.

We verify hypothesis
\begin{equation}\label{i29}
\mathbb{H}_0 \,\,:\,\, S(t) = 0, \quad t \in [0,1],
\end{equation}
versus simple alternatives
\begin{equation}\label{i30}
\mathbb{H}_\varepsilon  \,\,:\,\, S(t) = S_\varepsilon(t), \quad t \in [0,1], \quad \varepsilon > 0.
\end{equation}
It is clear that the asymptotic of type II error probabilities for arbitrary family of sets of alternatives $\Psi_\varepsilon \subset \mathbb{L}_2(0,1)$, $\varepsilon > 0$, is unambiguously described such a setup. We can always extract from the sets $\Psi_\varepsilon$, $\varepsilon > 0$, a family of simple alternatives  $S_\varepsilon$, $\varepsilon > 0$,  having worst order of asymptotic of type II error probabilities and to explore this asymptotic.

Paper is organized as follows. In section \ref{sec2}, we provide definitions of uniform consistency  in zones of large and moderate deviation probabilities. In section \ref{sec3} the asymptotics of type I and type II error probabilities in large and moderate deviations zones are pointed out if test statistics are linear combinations of squared estimates of Fourier coefficients of signal. In section \ref{sec4} we introduce the notion of maxisets for zones of large  deviation probabilities. Maxisets for test statistics are the largest convex, center-symmetric sets such that, if, from such a maxiset, we delete $\mathbb{L}_2$-balls, having centers at zero and having a given rates of convergence of radii to zero  as the noise power tends to zero, we get uniformly consistent family of sets of alternatives in large deviation sense. In section \ref{sec5} we establish that Besov bodies $\mathbb{B}^s_{2\infty}(P_0)$, $P_0 > 0$ are maxisets. Parameter $s$ depends on rates of convergence $\mathbb{L}_2$-norms $\|S_\varepsilon\|$ to zero. We  show that any family of simple alternatives $S_\varepsilon$, $\varepsilon > 0$, having given rates of convergence of $\mathbb{L}_2$-norms $\|S_\varepsilon\|$ to zero, is consistent in terms of large deviation  probabilities, if and only if, each function $S_\varepsilon$ admits representation as sum of two functions: $S_{1\varepsilon} \in \mathbb{B}^s_{2\infty}(P_0) $ with the same rates of convergence of $\|S_{1\varepsilon}\|$ to zero as $\|S_\varepsilon\|$ and orthogonal function $S_{2\varepsilon}$.  The functions
$S_{1\varepsilon}$ are smoother than $S_{2\varepsilon}$ and functions $S_{2\varepsilon}$  are more rapidly oscillating. In section \ref{sec6} we establish similar results if test statistics are $\mathbb{L}_2$-norms of kernel estimators of signal. In section \ref{sec7} proof of Theorems is provided.

We use letters $c$ and $C$ as a generic notation for positive constants. Denote ${\bf 1}_{\{A\}}$ the
indicator of an event $A$.  Denote $[a]$ whole part of real number $a$. For any two sequences of positive real numbers $a_n$ and $b_n$, $a_n = O(b_n)$ and $a_n \asymp b_n$ imply respectively $a_n < C b_n$ and $ca_n \le b_n \le Ca_n$ for all $n$ and $a_n = o(b_n)$, $a_n <<b_n$ implies $a_n/b_n \to 0$ as $n \to \infty$.  We also use notation   $a_n =  \Omega(b_n)$  which means  $b_n \le C a_n$ for all $n$. For any complex number $z$ denote $\bar z$ complex conjugate number.

Define function
$$ \Phi(x) = \frac{1}{\sqrt{2\pi}}\,\int_{-\infty}^x\,\exp\{-t^2/2\}\, dt, \quad x \in \mathbb{R}^1,
$$
of standard normal distribution.

Let $\phi_j$, $1 \le j < \infty$, be orthonormal system of functions onto $\mathbb{L}_2(0,1)$. For each $P_0 > 0$ define set
\begin{equation}\label{vv}
\mathbb{\bar B}^s_{2\infty}(P_0) = \Bigl\{S : S = \sum_{j=1}^\infty\theta_j\phi_j,\,\,\,  \sup_{\lambda>0} \lambda^{2s} \sum_{j>\lambda} \theta_j^2 \le P_0,\,\, \theta_j \in \mathbb{R}^1 \Bigr\}.
\end{equation}
With some conditions on the basis  $\phi_j$, $1 \le j < \infty,$ functional space
$$
\bar{\mathbb{ B}}^s_{2\infty} = \Bigl\{ S : S = \sum_{j=1}^\infty\theta_j\phi_j,\,\,\,  \sup_{\lambda>0} \lambda^{2s} \sum_{j>\lambda}\, \theta_j^2 < \infty,\,\, \theta_j \in \mathbb{R}^1 \Bigr\}
$$
is Besov space $\mathbb{B}^s_{2\infty}$ (see   \cite{rio}). In particular, this holds if
 $\phi_j$, $1 \le j < \infty$ is trigonometric basis.

If $\phi_j(t) = \exp\{2\pi i j x\}$, $x\in (0,1)$, $j = 0, \pm 1, \ldots$,  denote
$$
\mathbb{ B}^s_{2\infty}(P_0) = \Bigl\{S\,:\, f = \sum_{j=-\infty}^\infty \theta_j\phi_j,\,\,\,  \sup_{\lambda>0} \lambda^{2s} \sum_{|j| >\lambda} |\theta_j|^2 \le P_0 \Bigr\}.
$$
Since functions $\phi_j$ are complex here, then $\theta_j$ are complex numbers as well, and $\theta_j =  \bar\theta_{-j}$ for all $-\infty < j < \infty$.

Denote $\mathbb{ B}^s_{2\infty}$-- Banach space generated by balls $\mathbb{ B}^s_{2\infty}(P_0)$, $P_0 > 0$.
\section{Consistency in large and moderate deviation zone \label{sec2} }
For any test $L_\varepsilon$, $\varepsilon > 0$, denote $\alpha(L_\varepsilon) = \mathbf{E}_0(L_\varepsilon)$ -- its type I error probability and $\beta(L_\varepsilon,S)= \mathbf{E}_S(1 -L_\varepsilon)$ -- its type II error  probability for alternative  $S \in \mathbb{L}_2(0,1)$.

Let $r_\varepsilon \to \infty$ as $\varepsilon \to 0$. If  we have
\begin{equation}\label{us3}
|\log \beta(L_\epsilon,S_\varepsilon)|= \Omega(r_\varepsilon^2)
  \end{equation}
   for all tests $L_\varepsilon$, $\alpha(L_\varepsilon) = \alpha_\varepsilon$, such that   $r_\varepsilon^{-2}\,\log\alpha(L_\varepsilon) \to  0$ as $\varepsilon \to  0$, then we say that family of alternatives  $S_\varepsilon$, $\varepsilon > 0$,  is $r_\varepsilon^2$-consistent from large deviation viewpoint ($r_\varepsilon- LD$ -consistent).

    Note that any $r_\varepsilon^2- LD$- consistent family of alternatives $S_\varepsilon$ is consistent, while the reverse is not obligatory. If (\ref{us3}) does not hold, we say that family of alternatives is  $r_\varepsilon^2  - LD$  inconsistent.
\section{Asymptotic of type I and type II error probabilities of quadratic test statistics \label{sec3}}
Using orthonormal system of functions   $\phi_j$, $1 \le j < \infty$, we can rewrite stochastic differential equation  (\ref{vv1}) in terms of sequence model (see \cite{ib81})
\begin{equation}\label{q2}
y_{\varepsilon j} = \theta_{ j} +\varepsilon \xi_j, \quad 1 \le j < \infty,
\end{equation}
where $$y_{\varepsilon j} = \int_0^1 \phi_j dY_\varepsilon(t), \quad \xi_j = \int_0^1\,\phi_j\,dw(t) \quad \mbox{ and}  \quad \theta_j = \int_0^1 S\,\phi_j\,dt.$$

Denote $\yb_\varepsilon =  \{y_{\varepsilon j}\}_{j=1}^\infty$ and $\thetab = \{\theta_j\}_{j=1}^\infty$.

For alternatives $S_\varepsilon$ we put $\theta_{\varepsilon j} = \int_0^1 S_\varepsilon\,\phi_j\,dt$.

  We can consider $\thetab$ as a vector in Hilbert space  $\mathbb{H}$ with the norm $\|\thetab\| = \Bigl(\sum_{j=1}^\infty \theta_j^2\Bigr)^{1/2}$. In what follows, we shall implement the same notation of norm  $\| \cdot \|$ for the space $\mathbb{L}_2$ and for $\mathbb{H}$. The sense of this notation will be always clear from context.

We explore the test statistics of Neyman type tests
\begin{equation*}
 T_\varepsilon(Y_\varepsilon) = \varepsilon^{-2}\sum_{j=1}^\infty \varkappa_{\varepsilon j}^2 y_{\varepsilon j}^2  - \rho^2_\varepsilon,
\end{equation*}
where $\rho_\varepsilon^2 = \sum_{j=1}^\infty \varkappa_{\varepsilon j}^2$.

Suppose that coefficients $\varkappa_{\varepsilon j}^2$ satisfy the following assumptions.

\noindent{\bf A1.} For any $\varepsilon>0$ sequence $\varkappa^2_{\varepsilon j}$ is decreasing.

\noindent{\bf A2.} There holds
\begin{equation}\label{q5}
 C_1 < A_\varepsilon = \,\varepsilon^{-4}\,\sum_{j=1}^\infty \varkappa_{\varepsilon j}^4 < C_2.
 \end{equation}

  Denote $\varkappa_\varepsilon^2=\varkappa^2_{\varepsilon k_\varepsilon}$, where $k_\varepsilon = \sup\Bigl\{k: \sum_{j < k} \varkappa^2_{\varepsilon j} \le \frac{1}{2} \rho_\varepsilon^2 \Bigr\}$.

\noindent{\bf A3.}   There are $C_1$  and $\lambda >1$, such that for any $\delta > 0$ and for each $\varepsilon>0$ and all $k_{\varepsilon(\delta)} = [(1+\delta)k_\varepsilon]$ we have
 $
\varkappa^2_{\varepsilon,k_{\varepsilon(\delta)}} < C_1(1 +\delta)^{-\lambda}\varkappa_\varepsilon^2.
$

\noindent{\bf A4.} We have $\varkappa_{\varepsilon 1}^2  \asymp \varkappa_\varepsilon^2$ as $\varepsilon \to 0$. For any $c>1$  there is $C$ such that $\varkappa_{\varepsilon,[ck_\varepsilon]}^2 \ge C\varkappa_\varepsilon^2$ for all $\varepsilon > 0$.

 Define test
\begin{equation*}
L_\varepsilon= \mathbf{1}_{( A^{-1/2}_\varepsilon T_\varepsilon(Y_\varepsilon) > x_{\alpha_\varepsilon})}.
\end{equation*}
Denote $D_\varepsilon(S_\varepsilon) = \varepsilon^{-4} \sum_{j=1}^\infty \varkappa_{\varepsilon j}^2 \,\theta_{\varepsilon j}^2 = \varepsilon^{-4} \left(T_\varepsilon(S_\varepsilon) + \varepsilon^2 \rho^2_\varepsilon\right)$. Denote
\begin{equation*}
B_\varepsilon(S_\varepsilon) = \frac{D_\varepsilon(S_\varepsilon)}{(2A_\varepsilon)^{1/2}},
\end{equation*}

\begin{thm}\label{th1} Assume A1 -A4.

 Let $x_{\alpha_\varepsilon} \to \infty $ and let $x_{\alpha_\varepsilon} = o(k_{\varepsilon}^{1/6})$ as $\varepsilon \to 0$.

Then we have
\begin{equation}\label{us1}
 \alpha(L_\varepsilon) = (1 - \Phi(x_{\alpha_\varepsilon}))(1 + o(1)).
\end{equation}

Let $\alpha_\varepsilon \to 0 $,  $x_{\alpha_\varepsilon} = O(B_\varepsilon(S_\varepsilon))$ and let
\begin{equation}\label{eth2}
D_\varepsilon(S_\varepsilon)  \to \infty, \quad D_\varepsilon(S_\varepsilon) = o(k_\varepsilon^{1/6}), \quad B_\varepsilon(S_\varepsilon) -  x_{\alpha_\varepsilon} = \Omega(B_\varepsilon(S_\varepsilon))
\end{equation}
 as $\varepsilon \to 0$.

 Then we have
 \begin{equation}\label{us1n}
\beta(L_\varepsilon,S_\varepsilon) = \Phi (x_{\alpha_\varepsilon}- B_\varepsilon(S_\varepsilon))(1 + o(1)).
 \end{equation}

Let $x_{\alpha_\varepsilon} \to \infty $ and let $x_{\alpha_\varepsilon} = o(k_{\varepsilon}^{1/2})$ as $\varepsilon \to 0$.

Then we have
\begin{equation}\label{us2}
\sqrt{2|\log\alpha(L_\varepsilon)|} =x_{\alpha_\varepsilon}(1 +o(1)), \quad  \alpha_\varepsilon \to  0
\end{equation}
as $\varepsilon \to  0$.

Let $\alpha_\varepsilon \to 0 $,  $x_{\alpha_\varepsilon} = O(B_\varepsilon(S_\varepsilon))$  and let
\begin{equation} \label{eth3}
D_\varepsilon(S_\varepsilon)  \to \infty,\quad D_\varepsilon(S_\varepsilon) = o(k_\varepsilon^{1/2}), \quad B_\varepsilon(S_\varepsilon) -  x_{\alpha_\varepsilon} = \Omega(B_\varepsilon(S_\varepsilon))
\end{equation}
as $\varepsilon \to 0$.

Then we have
\begin{equation} \label{us2n}
\sqrt{2|\log\beta(L_\varepsilon,S_\varepsilon)|} = ( B_\epsilon(S_\varepsilon) - x_{\alpha_\varepsilon})(1 + o(1)),
\end{equation}
Let $\alpha_\varepsilon \to 0 $,  $x_{\alpha_\varepsilon} = O(B_\varepsilon(S_\varepsilon))$ and let
\begin{equation*}
r_\varepsilon \to \infty, \quad \frac{B_\varepsilon(S_\varepsilon)}{r_\varepsilon} \to \infty, \quad B_\varepsilon(S_\varepsilon) -  x_{\alpha_\varepsilon} = \Omega(B_\varepsilon(S_\varepsilon))
 \end{equation*}
 as $\varepsilon \to 0$.

 Then we have  $r_\varepsilon^{-2} |\log \beta(L_\varepsilon, S_\varepsilon)|  \to \infty$ as $\varepsilon \to 0$.\end{thm}

Denote
\begin{equation*}
\tau_\varepsilon^2= \varepsilon^4 D_\varepsilon(S_\varepsilon)=\sum_{j=1}^\infty \varkappa_{\varepsilon j}^2 \theta_{\varepsilon j}^2 = T_\varepsilon(S_\varepsilon) + \varepsilon^2\rho^2_\varepsilon.
\end{equation*}

\begin{thm}\label{th2} Assume $A1-A4$.

Let $x_{\alpha_\varepsilon} \to \infty $ and let  $k_\varepsilon = \Omega(x^2_{\alpha_\varepsilon})$  as $\varepsilon \to 0$. Then we have
\begin{equation}\label{eth20}
\log P_{0} (T_\varepsilon(Y_\varepsilon) > x_{\alpha_\varepsilon}) \le - \frac{1}{2} \varkappa_{\varepsilon 1}^{-2} x_{\alpha_\varepsilon} (1 + o(1)).
\end{equation}
Let $k_\varepsilon = \Omega(D_\varepsilon^2(S_\varepsilon))$  and $\varepsilon(\tau_\varepsilon - x_{\varepsilon}^{1/2}) \to \infty$ as $\varepsilon \to 0$ additionally. Then we have
\begin{equation}\label{eth21}
\log P_{S_\varepsilon} (T_\varepsilon (Y_\varepsilon) < x_{\alpha_\varepsilon}) \le - \frac{1}{2} \varkappa_{\varepsilon 1}^{-2} (\tau_{\varepsilon}-x^{1/2}_{\alpha_\varepsilon})^2(1 +o(1)).
\end{equation}
The equality in (\ref{eth21}) is attained for $S_{\varepsilon}= \{\theta_{\varepsilon j}\}_{j=1}^\infty$, where $\theta_{\varepsilon 1} = \varkappa_{\varepsilon 1}^{-1} \tau_\varepsilon$ and  $\theta_{\varepsilon j} =0$ for $j >1$.
\end{thm}
In Theorems \ref{th1} and \ref{th2} coefficients $\varkappa^2_{\varepsilon j}$ have different normalization in comparison with coefficients $\kappa^2_{\varepsilon j}$ in \cite{erm08}. At the same time  Theorems \ref{th1} and \ref{th2} are obtained by simple modification of proofs of Lemmas 2 and 3 in \cite{erm08}. To use notation $\varkappa^2_{\varepsilon j}$ instead of $\kappa^2_{\varepsilon j}$ in the proof of Lemmas 2 and 3 in \cite{erm08}, we should put
\begin{equation}\label{norm}
\kappa^2_{\varepsilon j} = \frac{\sum_{j=1}^\infty \varkappa^2 \theta^2_{\varepsilon j}}{\sum_{j=1}^\infty \kappa^4_{\varepsilon j}} \,\varkappa^2_{\varepsilon j}, \quad 1 \le j < \infty,
\end{equation}
Normalization (\ref{q5}) was implemented in Theorems on asymptotic normality of test statistics $T_\varepsilon(Y_\varepsilon)$ in \cite{er90, erm21}.
\section{Maxisets in large deviation zone \label{sec4} }
Maxiset definition in the zone of large deviation probabilities is akin to  the definition in  \cite{erm21} introduced for exploration of  uniform consistency in   problems of nonparametric hypothesis testing. The only difference is replacement of notion of consistency with notion of  $\varepsilon^{-2\omega}- LD$ consistency, $0 < \omega \le 1$.
 Just like in \cite{erm21}  we explore maxisets for alternatives  $S_\varepsilon$, approaching with hypothesis in  $\mathbb{L}_2$-- norm with the rate of convergence $\varepsilon^{2r}$, $ 0 < r < 1/2$, as $\varepsilon \to 0$.

Let $\Xi$, $\Xi \subset \mathbb{L}_2(0,1)$, be Banach space with the norm  $\|\cdot\|_\Xi$. Denote $ U=\{f:\, \|S\|_\Xi \le 1,\, S \in \Xi\}$ ball  in $\Xi$. In what follows, we suppose that set $U$ is compact in $\mathbb{L}_2(0,1)$ ( see \cite{ib77, erm21}).

Define subspaces $\Pi_k$, $1 \le k < \infty$, using induction.

Denote $d_1= \max\{\|S\|,\, S \in U\}$ and define function $e_1 \in U$ such that  $\|e_1\|= d_1.$ Denote $\Pi_1 \subset  \mathbb{L}_2(0,1)$ linear subspace generated by function $e_1$.

For $i=2,3,\ldots$ denote
$d_i = \max\{\rho(S,\Pi_{i-1}), S \in U \}$, where $\rho(S,\Pi_{i-1})=\min\{\|S-g\|, g \in \Pi_{i-1} \}$. Define function $e_i$, $e_i \in U$, such that $\rho(e_i,\Pi_{i-1}) = d_i$.
Denote $\Pi_i$ linear  space generated by functions $e_1,\ldots,e_i$.

For any function $S \in \mathbb{L}_2(0,1)$ denote
$S_{\Pi_i}$ projection of function  $S$ onto subspace $\Pi_i$ and put $\tilde S_i = S - S_{\Pi_i}$.

Set $U$ is called maxiset and functional space  $\Xi$ is called maxispace, if the following two statements hold:
\vskip 0.2cm
{\sl i.} any sequence of alternatives $S_\varepsilon \in U$, $ \|S_{\varepsilon}\| \asymp \varepsilon^{2r}$ as $\varepsilon \to 0$,  is $\varepsilon^{-2\omega}- LD$ consistent.
\vskip 0.2cm
{\sl ii.}  for any  $S \in \mathbb{L}_2(0,1)$, $S \notin \Xi$,  there is sequences $i_n$ and $j_{i_n}$, $i_n \to \infty$, $j_{i_n} \to \infty$ as $n \to \infty$   such that $c j_{i_n}^{-r}<\| \tilde S_{i_n}\| < C j_{i_n}^{-r}$ and subsequence of alternatives $\tilde S_{i_n}$ is $j_{i_n}^{\omega}- LD$ inconsistent for subsequence of signals $\tilde S_{i_n}$, observed in Gaussian white noise having the power $j_{i_n}^{-1/2}$.
\section{Necessary and sufficient conditions of $\varepsilon^{-2\omega}- LD$ - consistency of families of simple alternatives $S_\varepsilon$ converging to hypothesis in $\mathbb{L}_2$ \label{sec5}}
We explore necessary and sufficient conditions of  $\varepsilon^{-2\omega}- LD$ - consistency of families of simple alternatives  $S_\varepsilon$, $\varepsilon>0$,  such that $\|S_\varepsilon\| \asymp \varepsilon^{2r}$ as $\varepsilon \to 0$. Here  $0< r <1/2$ and $0 < 2\omega < 1 - 2r$. Thus $r_\varepsilon \asymp \varepsilon^{-\omega}$ as $\varepsilon \to 0$. We put $k_\varepsilon \asymp \varepsilon^{-4+ 8r +4\omega}$.

For this setup, A1-A4 implies
\begin{equation}\label{u1}
\kappa_\varepsilon^2  \asymp \varepsilon^{4-4r-4\omega}, \quad \bar A_\varepsilon \doteq \varepsilon^{-4}\,\sum_{j=1}^\infty \kappa_{\varepsilon j}^4\asymp \varepsilon^{-4\omega}
\end{equation}
as $\varepsilon \to 0$.

Theorems of section and subsequent proofs are akin to the statements and the proofs of Theorems  4.1 and 4.4-4.10 in \cite{erm21} with the above orders $\kappa_\varepsilon^2, \bar A_\varepsilon$ and $k_\varepsilon$. For the proofs it suffices to substitute the above orders in the proofs of corresponding Theorems in  \cite{erm21}. This will be demonstrated in the proof of the property {\it ii.} of maxiset in Theorem \ref{tq1}. All other proofs will be omitted.

 \subsection{Analytical form of necessary and sufficient conditions of $\varepsilon^{-2\omega}- LD$-consistency \label{s4.2}}
 Results are provided in terms of Fourier coefficients  of functions $S= S_\varepsilon = \sum_{j=1}^\infty \theta_{\varepsilon j} \phi_j$.
\begin{thm}\label{tq3} Assume {\rm A1-A4}. Family of alternatives $S_\varepsilon$, $\|S_\varepsilon\| \asymp \varepsilon^{2r}$, is  $\varepsilon^{-2\omega}- LD$ consistent, if and only if, there are $c_1$, $c_2$ and $\varepsilon_0>0$   such that there holds
\begin{equation}\label{con2}
\sum_{|j| < c_2k_\varepsilon} |\theta_{\varepsilon j}|^2 > c_1 \varepsilon^{4r}
\end{equation}
for all $\varepsilon < \varepsilon_0$.
\end{thm}
Versions of Theorem \ref{tq3} and Theorem \ref{tq6}, given bellow, hold also for test statistics based on $\mathbb{L}_2$-- norm of kernel estimators.  In this setup index  $j$ have positive and negative values and coefficients  $\theta_{\varepsilon j}$ may be complex numbers. By this reason, we write $|j|$ instead of $j$ and $|\theta_{\varepsilon j}|$ instead of $\theta_{\varepsilon j}$ in (\ref{con2}) and (\ref{con19}).
 \subsection{Maxisets. Qualitative structure of consistent families of alternatives. }

Denote $s = \frac{r}{2 -4r-2\omega}$. Then $r = \frac{(2-2\omega)s}{1 + 4s}$.
\begin{thm}\label{tq1} Assume {\rm A1-A4}. Then the balls  $\mathbb{\bar B}^s_{2\infty}(P_0)$, $P_0 > 0$, are maxisets for test statistics  $T_\varepsilon(Y_\varepsilon)$. \end{thm}
\begin{thm}\label{tq7} Assume {\rm A1-A4}. Then family of alternatives $S_\varepsilon$, $ \|S_\varepsilon\| \asymp \varepsilon^{2r}
$, is $\varepsilon^{-2\omega}- LD$ consistent, if and only if,  there is maxiset
$\mathbb{\bar B}^s_{2\infty}(P_0)$, $P_0 > 0$, and family of functions   $S_{1\varepsilon} \in
\mathbb{\bar B}^s_{2\infty}(P_0)$, $\|S_{1\varepsilon}\| \asymp \epsilon^{2r}$, such that $S_{1\varepsilon}$ is orthogonal
  $S_{\varepsilon} - S_{1\varepsilon}$, that is, we have
\begin{equation}
\label{ma1}
 \| S_\varepsilon\|^2 =  \| S_{1\varepsilon}\|^2  + \|S_\varepsilon -
S_{1\varepsilon}\|^2.
  \end{equation}
\end{thm}
\begin{thm}\label{tq11} Assume {\rm A1-A4}. L Then, for any $\delta > 0$, for any $\varepsilon^{-2\omega}- LD$ - consistent family of alternatives  $S_\varepsilon$,   $ \|S_\varepsilon\| \asymp \varepsilon^{2r}
$, there is maxiset  $\mathbb{\bar B}^s_{2\infty}(P_0)$, $P_0 > 0$, and family of functions $S_{1\varepsilon}$, $ \|S_{1\varepsilon}\| \asymp \varepsilon^{2r}
$,    $S_{1\varepsilon}  \in \mathbb{B}^s_{2\infty}(P_0)$,  such that there hold:
\vskip 0.25cm
$S_{1\varepsilon}$ is orthogonal  $S_\varepsilon - S_{1\varepsilon}$,
\vskip 0.25cm
  for any tests $L_\varepsilon$, satisfying (\ref{eth3}),
   there is $\varepsilon_0 =\varepsilon_0(\delta) > 0$, such that, for $\varepsilon < \varepsilon_0$, there hold
\begin{equation}
\label{uuu}
|\log \beta(L_\varepsilon,S_\varepsilon) - \log \beta(L_\varepsilon,S_{1\varepsilon})| \le \delta |\log \beta(L_\varepsilon,S_\varepsilon)|
\end{equation}
and
\begin{equation}
\label{uu1}
|\log \beta(L_\varepsilon,S_\varepsilon - S_{1\varepsilon}) | \le \delta |\log \beta(L_\varepsilon,S_\varepsilon)|.
\end{equation}
\end{thm}
\subsection{Interaction of $\varepsilon^{-2\omega}- LD$ consistent and $\varepsilon^{-2\omega}- LD$ inconsistent alternatives. Purely $\varepsilon^{-2\omega}- LD$ consistent families of alternatives}
We say that $\varepsilon^{-2\omega}- LD$- consistent family of alternatives   $S_\varepsilon$, $\|S_\varepsilon\| \asymp \varepsilon^{2r}$, is {\sl purely $\varepsilon^{-2\omega}- LD$--consistent}, if there is no inconsistent subsequence of alternatives $S_{1\varepsilon_i}$, $\varepsilon_i \to 0$ as $i \to \infty$, such that $S_{1\varepsilon_i}$ is orthogonal  $S_{\varepsilon_i}- S_{1\varepsilon_i}$ and $\|S_{1\varepsilon_i}\| > c_1\varepsilon_i^{2r}$.
\begin{thm}  \label{tq5}  Assume {\rm A1-A4}. Let family of alternatives $S_\varepsilon$, $\|S_\varepsilon\| \asymp \varepsilon^{2r}$, be $\varepsilon^{-2\omega}- LD$ consistent. Then for any  $\varepsilon^{-2\omega}- LD$ inconsistent family of alternatives $S_{1\varepsilon}$, $\|S_{1\varepsilon}\| \asymp \varepsilon^{2r}$,   there  holds
\begin{equation*}
\lim_{\varepsilon \to 0} \frac{\log\beta(L_\varepsilon,S_\varepsilon) - \log\beta(L_\varepsilon,S_\varepsilon + S_{1\varepsilon})}{\log\beta(L_\varepsilon,S_\varepsilon)} = 0.
\end{equation*}
\end{thm}
\begin{thm}\label{tq6} Assume {\rm A1-A4}. Family of alternatives $S_\varepsilon$, $\|S_\varepsilon\| \asymp \varepsilon^{2r}$, is purely $\varepsilon^{-2\omega}- LD$--consistent, if and only if, for any $\delta >0$,  there is such a constant $C_1= C_1(\delta)$, that there holds
\begin{equation}\label{con19}
\sum_{|j| > C_1k_\varepsilon} |\theta_{\varepsilon j}|^2 \le \delta \varepsilon^{4r}
\end{equation}
for all $\varepsilon< \varepsilon_0(\delta)$.
\end{thm}
\begin{thm}\label{tq12} Assume {\rm A1-A4}. Then family of alternatives $S_\varepsilon$, $ \|S_\varepsilon\| \asymp \varepsilon^{2r}
$, is purely $\varepsilon^{-2\omega}- LD$--consistent, if and only if, for any $\delta > 0$  there  is maxiset  $ \mathbb{\bar B}^s_{2\infty}(P_0)$ and family of functions $S_{1\varepsilon} \in \mathbb{\bar B}^s_{2\infty}(P_0)$   such that  $\|S_\varepsilon - S_{1\varepsilon}\| \le \delta \varepsilon^{2r}$ for all $\varepsilon< \varepsilon_0(\delta)$.
\end{thm}
\begin{thm}\label{tq8} Assume {\rm A1-A4}. Then family of alternatives $S_\varepsilon$, $\|S_\varepsilon\| \asymp \varepsilon^{2r}$, is purely $\varepsilon^{-2\omega}- LD$--consistent, if and only if, for any $\varepsilon^{-2\omega}- LD$--inconsistent sequence of alternatives $S_{1\varepsilon_i}$,  $ \|S_{1\varepsilon_i}\| \asymp \varepsilon_i^{2r}$, $\varepsilon_i \to 0$  as $i \to \infty$, there holds
\begin{equation}\label{ma2}
\|S_{\varepsilon_i} + S_{1\varepsilon_i}\|^2 =  \|S_{\varepsilon_i} \|^2 + \| S_{1\varepsilon_i}\|^2 + o(\varepsilon_i^{4r}).
\end{equation}
\end{thm}
\begin{remark}\label{rem1}{\rm Let $\varkappa_{\varepsilon j}^2 > 0$ for $j \le l_\varepsilon$,  and let $\varkappa^2_{\varepsilon j} = 0$ for $j > l_\varepsilon$, where $l_\varepsilon \asymp \varepsilon^{-4+ 8r +4\omega}$ as $\varepsilon \to 0$.
Analysis of proof of Theorems show that Theorems \ref{th1}, \ref{th2}  and \ref{tq3} - \ref{tq8}   remains valid for this setup, if assumption   A4  we replace with
\vskip 0.25cm
\noindent{\bf A5.} For any $c$, $0 < c <1$, there is $c_1$  such that $\varkappa^2_{\varepsilon,[cl_\varepsilon]} \ge c_1 \varkappa^2_{\varepsilon 1}$ for all $\varepsilon > 0$.
\vskip 0.25cm
In all corresponding proofs we put $\varkappa^2_\varepsilon = \varkappa_{\varepsilon 1}^2$ and $k_\varepsilon = l_\varepsilon$.
Theorems   \ref{tq6} is correct with the following modifications. It suffices  suppose that $C_1(\epsilon) < 1$.
Proof of corresponding versions of Theorems\ref{th1}, \ref{th2}  and \ref{tq3} - \ref{tq8} is similar and is omitted.}
\end{remark}

\section{Tests based on kernel estimators \label{sec6}}
We explore the same problem of signal detection in Gaussian white noise (\ref{i29}), (\ref{i30}). Suppose additionally,  that signal  $S_\varepsilon$   belong to $\mathbb{L}_2^{per}(\mathbb{R}^1)$ the set of 1-periodic functions such that $S_\varepsilon(t) \in \mathbb{L}_2(0,1)$, $t \in [0,1)$. This allows to extend our model on real line $\mathbb{R}^1$ putting $w(t+j) = w(t)$ for all integer $j$ and $t \in [0,1)$ and  to write the forthcoming integrals over all real line (see \cite{erm21}).

Define kernel estimators
\begin{equation}\label{yy}
\hat{S}_\varepsilon(t) = \frac{1}{h_\varepsilon} \int_{-\infty}^{\infty} K\Bigl(\frac{t-u}{h_\varepsilon}\Bigr)\, d\,Y_\varepsilon(u), \quad t \in (0,1),
\end{equation}
where $h_\varepsilon> 0$, $h_\varepsilon \to 0$ as $\varepsilon \to 0$.

The kernel $K$ is bounded function such that the support of $K$ is contained in $[-1/2,1/2]$, $K(t) = K(-t)$ for $t \in \mathbb{R}^1$ and $\int_{-\infty}^\infty K(t)\,dt = 1$.

Denote $K_h(t) = \frac{1}{h} K\Bigl(\frac{t}{h}\Bigr)$, $t \in \mathbb{R}^1$ and $h >0$.

Define test statistics $$T_{1\varepsilon}(Y_\varepsilon) = T_{1\varepsilon h_\varepsilon}(Y_\varepsilon) =\varepsilon^{-2}h_\varepsilon^{1/2}\gamma^{-1} (\|\hat S_{\varepsilon}\|^2- \varepsilon^2 h_\varepsilon^{-1}\|K\|^2),$$  where
$$
\gamma^2 = 2 \int_{-\infty}^\infty \Bigl(\int_{-\infty}^\infty K(t-s)K(s) ds\Bigr)^2\,dt.
$$
Statistics $\|\tilde S_\varepsilon\|^2- \varepsilon^2 h_\varepsilon^{-1}\|K\|^2$ is estimator of functional value
$$
T_{\varepsilon}(S_\varepsilon)  =\int_0^1\Bigl(\frac{1}{h_\varepsilon}\int K\Bigl(\frac{t-s}{h_\varepsilon}\Bigr)S_\varepsilon(s)\, ds\Bigr)^2 dt.
$$
Denote $L_\varepsilon = \mathbf{1}_{\{T_{1\varepsilon}(Y_\varepsilon) > x_{\alpha_\varepsilon}\}}$
, where $x_{\alpha_\varepsilon}$ is defined by equation $\alpha(L_\varepsilon) = \alpha_\varepsilon$.

Problem admits interpretation in the framework of sequence model. Let we observe realization of random process
 $Y_\varepsilon(t)$, $t  \in [0,1]$.

 For $-\infty < j < \infty$, denote
$$
\hat K(jh) = \int_{-1}^1 \exp\{2\pi ijt\}\,K_h(t)\, dt,\quad h > 0,
$$
$$
y_{\varepsilon j} = \int_0^1 \exp\{2\pi ijt\}\, dY_\varepsilon(t),
\quad
\xi_j = \int_0^1 \exp\{2\pi ijt\}\, dw(t),
$$
$$
\theta_{\varepsilon j} = \int_0^1 \exp\{2\pi ijt\}\, S_\varepsilon(t)\, dt.
$$
In this notation we can define estimator in the following form
\begin{equation}\label{au1}
\hat \theta_{\varepsilon j} = \hat K(jh_\varepsilon)\, y_{\varepsilon j} =
\hat K(jh_\varepsilon)\, \theta_{\varepsilon j} +  \varepsilon\, \hat K(jh_\varepsilon)\, \xi_j, \quad -\infty < j < \infty.
\end{equation}
Then test statistics $T_\varepsilon$ admits the following representation:
\begin{equation}\label{au111}
T_\varepsilon(Y_\varepsilon) = \varepsilon^{-2} h_\varepsilon^{1/2}\gamma^{-1} \Bigl(\sum_{j=-\infty}^\infty |\hat \theta_{\varepsilon j}|^2  -  \varepsilon^2 \sum_{j=-\infty}^\infty |\hat K(jh_\varepsilon)|^2\Bigr).
\end{equation}
If we put $|\hat K(jh_\varepsilon)|^2 = \kappa^2_{\varepsilon j}$, we get, that test statistics  $T_\varepsilon(Y_\varepsilon)$ in this section and sections \ref{sec3} and \ref{sec5} are almost indistinguishable. Results  similar to section \ref{sec3} were obtained in \cite{erm11}. However for obtaining results similar to section   \ref{sec5} one needs their interpretation in terms of Fourier coefficients.

Main difference of setup of section   \ref{sec6} is the presence of heterogeneous Gaussian white noise in the model.
\begin{thm}\label{tk3} Let $\|S_\varepsilon\| \asymp \varepsilon^{2r}$  and $h_\varepsilon \asymp \varepsilon^{4-8r-4\omega}$, $ 0 < 2 \omega < 1 - 2r$. Then Theorems \ref{tq3}--\ref{tq8} are valid for this setup with the only difference that $ \mathbb{\bar B}^s_{2\infty}(P_0)$  is replaced with $\mathbb{B}^s_{2\infty}(P_0)$.
\end{thm}
Proof of Theorems is based on the following version of Theorem  \ref{th1} (see Theorems  2.1 and 2.2,
 \cite{erm11}).

\begin{thm}\label{tk2} Let $h_\varepsilon \to 0$ as $\varepsilon \to 0$.
Let $\alpha_\varepsilon \doteq \alpha(K_\varepsilon) = o(1)$.

If
\begin{equation*}
1 < < \sqrt{2|\log \alpha_\varepsilon|} < <   h_\varepsilon^{-1/6},
\end{equation*}
then $x_{\alpha_\varepsilon}$ is defined equation $\alpha_\varepsilon= (1 - \Phi(x_{\alpha_\varepsilon}))(1 + o(1))$ as $\varepsilon \to 0$.

If
\begin{equation}\label{metka}
\varepsilon^{-2} h_\varepsilon^{1/2} T_\varepsilon(S_\varepsilon) - \sqrt{|2\log \alpha_\epsilon|} > c \varepsilon^{-2} h_\varepsilon^{1/2} T_\varepsilon(S_\varepsilon)
\end{equation}
and
\begin{equation*}
\varepsilon^2 h_\varepsilon^{-1/2} < < T_\varepsilon(S_\varepsilon) < < \varepsilon^2
h_\varepsilon^{-2/3} = o(1),
\end{equation*}
as $\varepsilon \to 0$, then
 the following relation holds
\begin{equation}\label{33}
\beta(L_\varepsilon,S_\varepsilon)  = \Phi(x_{\alpha_\varepsilon} - \gamma^{-1} \varepsilon^{-2} h_\varepsilon^{1/2} T_{\varepsilon}(S_\varepsilon))(1 + o(1)).
\end{equation}
If
\begin{equation*}
1 < < \sqrt{2|\log \alpha_\varepsilon|} < <   h_\varepsilon^{-1/2},
\end{equation*}
then $x_{\alpha_\varepsilon} =  \sqrt{|2\log\alpha_\varepsilon|}(1 + o(1))$.

If
\begin{equation*}
\varepsilon^2 h_\varepsilon^{-1/2} < < T_\varepsilon(S_\varepsilon) < < \varepsilon^2  h_\varepsilon^{-1} = o(1),
\end{equation*}
and (\ref{metka}) holds then we have
\begin{equation}\label{33}
2\log\beta(L_\varepsilon,S_\varepsilon)  = -(x_{\alpha_\varepsilon} - \gamma^{-1} \varepsilon^{-2} h_\varepsilon^{1/2} T_{\varepsilon}(S_\varepsilon))^2(1 + o(1)).
\end{equation}
If $\frac{\varepsilon^2 x_{\alpha_\varepsilon}}{h_\varepsilon^{1/2} T_{\varepsilon}(S_\varepsilon)} \to 0 $ as $\varepsilon \to 0$, then
\begin{equation*}
\lim_{\varepsilon \to 0} (\log \alpha_\varepsilon)^{-1} \log \beta (L_\varepsilon,S_\varepsilon) = \infty.
\end{equation*}
\end{thm}
Note that
\begin{equation}\label{z33}
T_{\varepsilon}(S_\varepsilon)   = \sum_{j=-\infty}^\infty |\hat K(jh_\varepsilon)|^2 |\theta_{\varepsilon j}|^2.
\end{equation}
\section{Proof of Theorems \label{sec7}}
\subsection{Proof of Theorem  \ref{th2}}
Denote $z^2 = \kappa^{2}_{\varepsilon} x_{\alpha_\varepsilon}$.
Implementing Chebyshev inequality, we have
\begin{equation}\label{dth21}
\begin{split} &
\mathbf{P}_s(T_\varepsilon - \tau_\varepsilon^2 < z^2 - \tau_\varepsilon^2)\\&
\le \exp\{t(z^2 - \tau_\varepsilon^2)\} \mathbf{E}_0 \left[\exp\left\{-t\sum_{j=1}^\infty \kappa_j^2(y_j - \varepsilon^2) + 2t\sum_{j=1}^\infty\kappa_j^2y_j\theta_j\right\}\right]\,d\,y \\&=
\exp\{t(z^2 - \tau_\varepsilon^2)\} \\& \times \lim_{m\to \infty} (2\pi\varepsilon^2)^{-m/2}\int \left[\exp\left\{-\frac{1}{2}\sum_{j=1}^m\left(\frac{1 + 2t\kappa_j^2\varepsilon^2}{\varepsilon^2} y_j^2 \right)\right.\right.\\&\left.\left. + t \sum_{j=1}^m \kappa_j^2 \epsilon^2+ 2t\sum_{j=1}^m\kappa_j^2y_j\theta_j \pm \sum_{j=1}^m \frac{2t^2\kappa_j^4\theta_j^2\varepsilon^2}{1 + 2t\kappa_j^2\varepsilon^2}   \right\}\right]\,d\,y \\& = \exp\{t(z^2 - \tau_\varepsilon^2)\}\exp\left\{\frac{1}{2}\sum_{j=1}^\infty\log(1 + 2t\kappa_j^2\varepsilon^2)  + \sum_{j=1}^\infty \frac{2t^2\kappa_j^4\theta_j^2\varepsilon^2}{1 + 2t\kappa_j^2\varepsilon^2}  \right\}.
\end{split}
\end{equation}
Note  that
\begin{equation*}
\sum_{j=1}^\infty\log (1 + 2t\kappa_j^2\varepsilon^2) \asymp k_\varepsilon << \varepsilon^{-2}\tau_\varepsilon^2.
\end{equation*}
It is easy to see that, if $\sum_{j=1}^\infty \kappa_j^2 \theta_j^2 = \mbox{const}$, then supremum on  $\theta_j$ of formula  under the sign of the second exponent is reached when $\kappa^2 \theta_1^2 = \tau_\varepsilon^2$.

Thus the problem is reduced to maximization on  $t$
\begin{equation}\label{dth22}
\exp\Bigl\{t(z^2 - \tau_\varepsilon^2) + \frac{2t^2\kappa_\varepsilon^2\tau_\varepsilon^2\varepsilon^2}{1 + 2t\kappa_\varepsilon^2\varepsilon^2}\Bigl\}.
\end{equation}
By straightforward calculations, we get that maximum on  $t$ is attained if
\begin{equation}\label{dth23}
2t = \varepsilon^{-2}\kappa_\varepsilon^{-2}(\tau_\varepsilon z^{-1} - 1).
\end{equation}
Substituting (\ref{dth23}) to (\ref{dth22}), we get right hand-side of  (\ref{eth21}).
\subsection{Proof of Theorem \ref{th1}}
Test statistics $T_\varepsilon$ are sums of independent random variables. Therefore the standard reasoning of proof of Cramer Theorem are implemented \cite{os,petr}. This reasoning has been realized in proof of Lemma 4 in \cite{erm08}. If hypothesis is valid, then proof of Theorem  \ref{th1} completely coincides  with proof of Lemma 4 in \cite{erm08}.
If alternative is valid,  the difference is only in the evaluation of the residual terms arising in the proof of Lemma 4 in \cite{erm08}.

Such differences will arise first of all in the evaluation of the residual term in (3.51) in \cite{erm08}.

If  (\ref{eth3}) holds, using $\varkappa_\varepsilon^2 \asymp \varepsilon^2 k_\varepsilon^{-1/2}$, we have
\begin{equation}\label{dth11a}
\varepsilon^{-6}\,\sum_{j=1}^\infty \varkappa_{\varepsilon j}^4\, \theta_{\varepsilon j}^2 \asymp \varepsilon^{-6}\,\varkappa_\varepsilon^2\,\sum_{j=1}^\infty \varkappa_{\varepsilon j}^2\, \theta_{\varepsilon j}^2 \asymp k_\varepsilon^{-1/2} \,D_\varepsilon(S_\varepsilon),
\end{equation}
or, for $\kappa_{\epsilon j} $-- notation, we have
\begin{equation}\label{dth11}
\varepsilon^{-6}\sum_{j=1}^\infty \kappa_{\varepsilon j}^4\, \theta_{\varepsilon j}^2 \asymp \varepsilon^{-6}\,\kappa_\varepsilon^2\,\sum_{j=1}^\infty \kappa_{\varepsilon j}^2\, \theta_{\varepsilon j}^2 = o\left(\varepsilon^{-4}\,\kappa_\varepsilon^2\,\sum_{j=1}^\infty \kappa_{\varepsilon j}^4\right).
\end{equation}

Using (\ref{dth11}),  we can replace remainder terms in  (3.53), (3.55)--(3.57) in \cite{erm08} by the factor $(1 +o(1))$. This allows to replace (3.43)  (see Lemma 4 in \cite{erm08}) with (\ref{us2n}).

If (\ref{eth3}) holds, then we have
\begin{equation}\label{dth12a}
\begin{split} &
( B_\varepsilon(S_\varepsilon)-x_{\alpha_\varepsilon})\,\varepsilon^{-6}\,\sum_{j=1}^\infty \varkappa_{\varepsilon j}^4 \,\theta_{\varepsilon j}^2 \,\asymp \, \varepsilon^{-2}\varkappa_\varepsilon^2\left(\varepsilon^{-4}\,\sum_{j=1}^\infty \varkappa_{\varepsilon j}^2 \,\theta_{\varepsilon j}^2\right)^3 \\& \asymp \,
k_\varepsilon^{1/2}\,\left(\varepsilon^{-4}\,\sum_{j=1}^\infty \varkappa_{\varepsilon j}^2\, \theta_{\varepsilon j}^2\right)^3 = o(1),
\end{split}
\end{equation}
or, for $\kappa_{\epsilon j} $-- notation, we have
\begin{equation}\label{dth12}
\begin{split} &
\varepsilon^{-6}\sum_{j=1}^\infty \kappa_{\varepsilon j}^4\, \theta_{\varepsilon j}^2 \,\asymp \, \varepsilon^{-6}\,\kappa_\varepsilon^2\,\sum_{j=1}^\infty \kappa_{\varepsilon j}^2 \, \theta_{\varepsilon j}^2 \\& \asymp \varepsilon^{-6}\,\kappa_\varepsilon^2\,\sum_{j=1}^\infty \, \kappa_{\varepsilon j}^4\asymp \varepsilon^{-6}\,k_\varepsilon\,\kappa_\varepsilon^6\, = \,o(1).
\end{split}
\end{equation}
Using (\ref{dth12}),  we can remain all estimates of the residual terms in the proof of Lemma 4 in \cite{erm08} and get the validity of these estimates for proof of  (\ref{us1n}).

Proof of (\ref{us1}) and (\ref{us2}) for type I error probabilities is based on similar estimates. It suffices to put $x_{\alpha_\varepsilon} = \varepsilon^{-4}\,\sum_{j=1}^\infty \varkappa_{\varepsilon j}^2 \,\theta_{\varepsilon j}^2$ in above mentioned estimates.

\subsection{Proof of Theorem \ref{tq1}}
We verify only {\it ii.} in definition of maxiset.  Proof of  {\it i.} is akin to \cite{erm21} and
is omitted.

 Suppose the contrary. Then  $S = \sum_{j=1}^\infty \tau_{j}\,\phi_j  \notin \mathbb{\bar
 B}^s_{2\infty}$. This implies that there is sequence  $m_l$, $m_l \to \infty$ as $l \to \infty$,
 such that
\begin{equation}\label{u5}
m_l^{2s} \sum_{j=m_l}^\infty \tau_{j}^2 = C_l,
\end{equation}
where $C_l \to \infty$ as $l \to \infty$.

Define sequence $\etab_l = \{\eta_{lj}\}_{j=1}^\infty$ such that
$\eta_{lj} = 0$ if $j  < m_l$ and $\eta_{lj} = \tau_j$, if $j \ge m_l$.

We put $\tilde S_l = \sum_{j=1}^\infty \eta_{lj}\,\phi_j$.

For alternatives $\tilde S_l$ we define sequence $n_l$ such that
\begin{equation}\label{u5b}
\|\etab_l\|^2 \asymp \varepsilon_l^{4r}= n_l^{-2r}\asymp m_l^{-2s} C_l .
\end{equation}
Then
\begin{equation}\label{u7}
n_l \asymp C_l^{-1/(2r)} m_l^{s/r} \asymp C_l^{-1/(2r)} m_l^{\frac{1}{2 - 4r- 2 \omega}}.
\end{equation}
Hence we have
\begin{equation}\label{u10}
m_l \asymp C_l^{(1-2r-\omega)/r}n_l^{2-4r-2\omega}.
\end{equation}
By A4, (\ref{u10}) implies
\begin{equation}\label{u9}
\kappa^2_{n_l m_l} = o(\kappa^2_{n_l}).
\end{equation}
Using (\ref{u1}), A2 and (\ref{u9}), we get
\begin{equation}\label{u12}
\begin{split}&
A_{\varepsilon_l}(\etab_l) = n_l^2 \sum_{j=1}^{\infty} \kappa_{\varepsilon_lj}^2\eta_{jl}^2 \le
n_l^{2} \kappa^2_{\varepsilon_l m_l}\sum_{j=m_l}^{\infty} \theta_{n_lj}^2\\& \asymp
n_l^{2-2r}\kappa^2_{\varepsilon_l m_l} \asymp \varepsilon_l^{-4\omega}\kappa^2_{\varepsilon_l m_l
}\kappa^{-2}_{\varepsilon_l}) = o(\varepsilon_l^{-4\omega}) .
\end{split}
\end{equation}
By Theorem \ref{th1},  (\ref{u12}) implies $\varepsilon^{-2\omega}- LD$ --inconsistency of sequence
of alternatives $\tilde S_l$.
\subsection{Proof of version of Theorem \ref{tq1} for tests based on kernel estimators}
We verify only {\it ii.} Suppose the contrary.
Let $S = \sum_{j=-\infty}^\infty \tau_j \phi_j \notin \gamma U$ for all $\gamma >0$.  Then there is
sequence  $m_l$, $m_l \to \infty$ as $l \to \infty$, such that
\begin{equation}\label{bb5}
m_l^{2s} \sum_{|j|\ge m_l}^\infty |\tau_j|^2 = C_l,
\end{equation}
where $C_l \to \infty$ as $l \to \infty$.

It is easy to see (see \cite{erm21}) that we can define sequence $m_l$ such that
\begin{equation}\label{gqq}
 \sum_{m_l \le |j| \le 2m_{l}} |\tau_j|^2 \asymp \varepsilon_l^{4r} = n_l^{-2r}> \delta C_l
 m_l^{-2s},
\end{equation}
where $\delta$, $0< \delta <1/2$,  does not depend on $l$.

Define sequence $\etab_l = \{\eta_{lj}\}_{j=-\infty}^\infty$ such that
 $\eta_{lj} = \tau_j$,  $|j| \ge m_{l}$, and $\eta_{lj} = 0$ otherwise.

Denote
$$
\tilde S_l(x) =  \sum_{j=-\infty}^\infty \eta_{lj} \exp\{2\pi ijx\}.
$$
For alternatives $\tilde S_l(x)$ we define sequence $n_l$ such that $\|\tilde S_l(x)\|
\asymp n_l^{-r}$.

 Then
\begin{equation*}
n_l \asymp C_l^{-1/(2r)} m_l^{s/r}.
\end{equation*}
We have $|\hat K(\omega)| \le \hat K(0) = 1$ for all $\omega \in R^1$ and $|\hat K(\omega)| > c >
0$ for all $ |\omega| < b$. Therefore, if we put $h_l= h_{n_l} =2^{-1}b^{-1}m_l^{-1}$,
then, by (\ref{gqq}), there is $C > 0$  such that for all $h> 0$, we have
\begin{equation*}
T_{\varepsilon_l}(\tilde S_l,h_l) = \sum_{j=-\infty}^\infty |\hat K(jh_l)\,\eta_{lj}|^2 > C
\sum_{j=-\infty}^\infty |\hat K(jh)\,\eta_{lj}|^2 = C T_{\varepsilon_l}(\tilde S_l,h).
\end{equation*}
Therefore, we can put  $h = h_l$ in further reasoning.

Using (\ref{gqq}), we get
\begin{equation}\label{k101}
T_{\varepsilon_l}(\tilde S_l) =  \sum_{|j|>m_l}\, |\hat K(jh_l) \eta_{lj}|^2  \asymp
\sum_{j=m_l}^{2m_l} |\eta_{lj}|^2 \asymp n_l^{-2r}.
\end{equation}
If we put in (\ref{u7}),(\ref{u10}), $k_l = [h_{\varepsilon_l}^{-1}]$ and $m_l = k_l$,
then we get
\begin{equation}\label{k102}
h_{\varepsilon_l}^{1/2} \asymp C_l^{(2r-1+1)/(2r)}n_l^{2r-1+\omega}.
\end{equation}
Using (\ref{k101}) and (\ref{k102}), we get
\begin{equation*}
\varepsilon_l^{-2} T_{\varepsilon_l}(\tilde S_l)h_{\varepsilon_l}^{1/2} \asymp C_l^{-(1-2r)/2}.
\end{equation*}                                                                                                                                                                                           By Theorem \ref{tk2}, this implies inconsistency of sequence of alternatives  $\tilde S_l$.

\end{document}